\documentclass[graybox, envcountchap]{svmult}

\usepackage{makeidx}         
\usepackage{graphicx}        
\usepackage{multicol}        
\usepackage[bottom]{footmisc}

\usepackage{newtxtext}       %
\usepackage[varvw]{newtxmath}       
\usepackage{amsmath}
\usepackage{mathtools, cuted}
\usepackage{siunitx}
\usepackage{booktabs}
\usepackage{multicol}
\usepackage{multirow}
\usepackage{bm}
\usepackage{array}
\newcolumntype{P}[1]{>{\centering\arraybackslash}p{#1}}
\newcolumntype{M}[1]{>{\centering\arraybackslash}m{#1}}
\usepackage[colorlinks,urlcolor=black,linkcolor=black,citecolor=black]{hyperref}

\makeindex             

\begin{document}

\frontmatter



\mainmatter
%
%
%

%
%
%
%
%
%
%

\title{Machine Learning Based Optimization Workflow for Tuning Numerical Settings of Differential Equation Solvers for Boundary Value Problems}
\titlerunning{ML-based optimization for tuning BVP solver settings}
\author{Viny Saajan Victor, Manuel Ettm\"{u}ller, Andre Schmei{\ss}er, Heike Leitte and Simone Gramsch}
\authorrunning{V. S. Victor, M. Ettm\"{u}ller, A. Schmei{\ss}er, H. Leitte and S. Gramsch}
\institute{Viny Saajan Victor \at Fraunhofer ITWM, Fraunhofer-Platz 1, 67663 Kaiserslautern, Germany, \email{viny.saajan.victor@itwm.fraunhofer.de}
\and Manuel Ettm\"{u}ller \at Fraunhofer ITWM, Fraunhofer-Platz 1, 67663 Kaiserslautern, Germany \email{manuel.ettmueller@itwm.fraunhofer.de}
\and Andre Schmei{\ss}er \at Fraunhofer ITWM, Fraunhofer-Platz 1, 67663 Kaiserslautern, Germany \email{andre.schmeisser@itwm.fraunhofer.de}
\and Heike Leitte \at Technical University of Kaiserslautern, Gottlieb-Daimler-Straße 47, 67663 Kaiserslautern, Germany \email{leitte@cs.uni-kl.de} 
\and Simone Gramsch \at Fraunhofer ITWM, Fraunhofer-Platz 1, 67663 Kaiserslautern, Germany \email{simone.gramsch@itwm.fraunhofer.de}}
%
%
\maketitle

\abstract{ Several numerical differential equation solvers have been employed effectively over the years as an alternative to analytical solvers to quickly and conveniently solve differential equations.
One category of these is boundary value solvers, which are used to solve real-world problems formulated as differential equations with boundary conditions. These solvers require certain numerical settings to solve the differential equations that affect their solvability and performance. A systematic fine-tuning of these settings is required to obtain the desired solution and performance. Currently, these settings are either selected by trial and error or require domain expertise. In this paper, we propose a machine learning-based optimization workflow for fine-tuning the numerical settings to reduce the time and domain expertise required in the process. In the evaluation section, we discuss the scalability, stability, and reliability of the proposed workflow. We demonstrate our workflow on a numerical boundary value problem solver.}

\section{Introduction}
\label{sec:1}
Differential equations are one of the most effective mathematical tools for understanding and predicting the behavior of dynamical systems in nature, engineering, and society. The study of differential equations entails learning how to solve them and interpret the solutions. In the field of differential equations, a boundary value problem (BVP) is a type of differential equation with a set of additional constraints known as boundary conditions. The solution to a boundary value problem must satisfy the boundary conditions. Despite the fact that analytical methods for solving boundary value problems yield exact answers, they become difficult to apply to complex problems. Hence, several numerical methods have become popular, leveraging the development of computing capabilities. These methods provide approximate solutions that have sufficient accuracy for engineering purposes. The boundary value problem solvers developed based on the numerical methods usually consist of settings that describe method-specific properties such as error tolerance, mesh size, etc. These settings have an impact on the solvability and performance of the solvers, and the impact varies depending on the problem. In order to obtain qualitatively good solutions within acceptable solution times, fine-tuning these settings is required. However, these settings are now tuned manually by trial and error, which is ineffective due to the large parameter space of the settings and the complex interactions between them.

Machine learning (ML) techniques have achieved remarkable success in recent years due to their ability to identify patterns and structures in the given data. Since ML models can predict in real-time, they have been deployed in a variety of industries. In earlier work, we showed how ML-based quality optimization is effective in the technical textile industry~\cite{victor2022visual}. In this paper, we propose a two-stage optimization workflow based on ML to make the process of tuning the numerical settings of boundary value problem solvers more effective. In the first stage of the workflow, we design an ML pipeline that maps the numerical settings to solver performance. The solver performance metric includes solvability status and computational cost statistics of the solver, such as the number of evaluations and mesh points required by the solver and the corresponding residual error obtained from the solver. In the second stage of the workflow, we use the trained ML model to predict the influence of the numerical settings on the success and performance of the solver. We use this knowledge to propose a multi-criteria optimization strategy for fine-tuning the settings, thereby reducing the domain expertise required in the process.

\section{Foundations: The test bench for boundary value problem solvers}
\label{sec:2}

In collaborations, we often work together with domain experts in the technical textile industry such as mathematicians or process engineers. They want to model and simulate a variety of fiber formation processes such as polymer melt spinning~\cite{ettmuller2021flow} or the production of glass wool~\cite{arne2011fluid}. These processes are typically modeled as a system of differential equations that require a numerical solver. In the special case of ordinary differential equations, the bvp4c solver from MATLAB is frequently used when dealing with boundary value problems~\cite{kierzenka2001bvp}.
The bvp4c solver is a finite difference method based on a collocation method that provides a continuous solution that is fourth-order accurate in the interval of integration. The integration interval is thereby subdivided into smaller intervals using a mesh of points. The subdivision scheme is adaptive, so further points can be both added and removed. The boundary conditions and the collocation conditions imposed on all of the sub-intervals create a global system of algebraic equations. These equations are solved using the simplified Newton method to arrive at a numerical solution.
For the formulation of the ML problem statement and demonstration of the proposed approach, we used a self-implemented C++ version of the solver above, which uses a Newton-Armijo method instead of a simplified Newton method for solving the system of non-linear equations.

\begin{figure}
    \centering
    \includegraphics[scale=0.8]{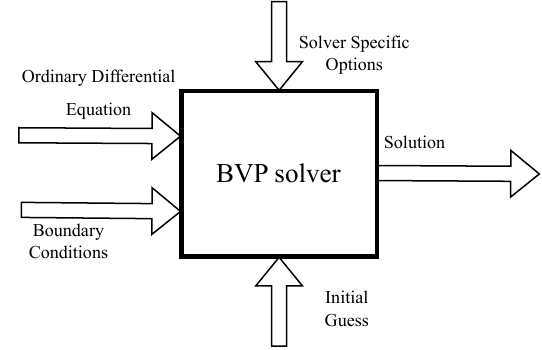}
    \caption{Input and output parameters of BVP solver.}
    \label{fig:bvp4c}
\end{figure}

The solvability status and the computational resources required by the BVP solver are mainly affected by the type of differential equations, boundary conditions, the initial guess of the solution, and the solver-specific options as shown in Figure~\ref{fig:bvp4c}. In our workflow, we keep the initial guess for the solver constant, as we are interested in analyzing the effect of solver-specific options. We also want to analyze the
correlation between the optimal choice of numerical settings and the complexity of the boundary value problem. For this, we selected ten different reference problems for our analysis. Out of all the test problems defined in~\cite{soetaert2010package}, we have chosen four linear (1, 3, 4, and 7) and six non-linear (19, 20, 22, 23, 24, and 33) problems, since we expect to gain more insights from the more complex differential equations.

\section{Machine learning based optimization workflow for tuning the BVP solver settings}
\label{sec:3}
In this section, we elaborate on the proposed workflow for tuning the numerical settings of the BVP solvers. The work consists of the machine learning process and the optimization process. The machine-learning process is used to map the numerical settings of the BVP solvers to their performance. This pipeline comprises 
a binary classifier and a multiple-output regression model. The binary classifier is used to predict the solvability status of the boundary value problem based on the numerical settings. The settings predicted to be solvable by the classifier are then fed into the multiple-output regression model which forecasts solver performance. Both machine learning models rely on the training data collected for the reference boundary value problems. The best models are chosen based on the accuracy obtained from the testing data. The chosen models are used to solve the multi-objective optimization problem to obtain the desired combination of numerical settings that optimize the quality criteria. We discuss the details of feature selection, data collection, model selection, training, and testing of ML models in the following subsections.

\begin{figure}
    \centering
    \includegraphics[width=\linewidth]{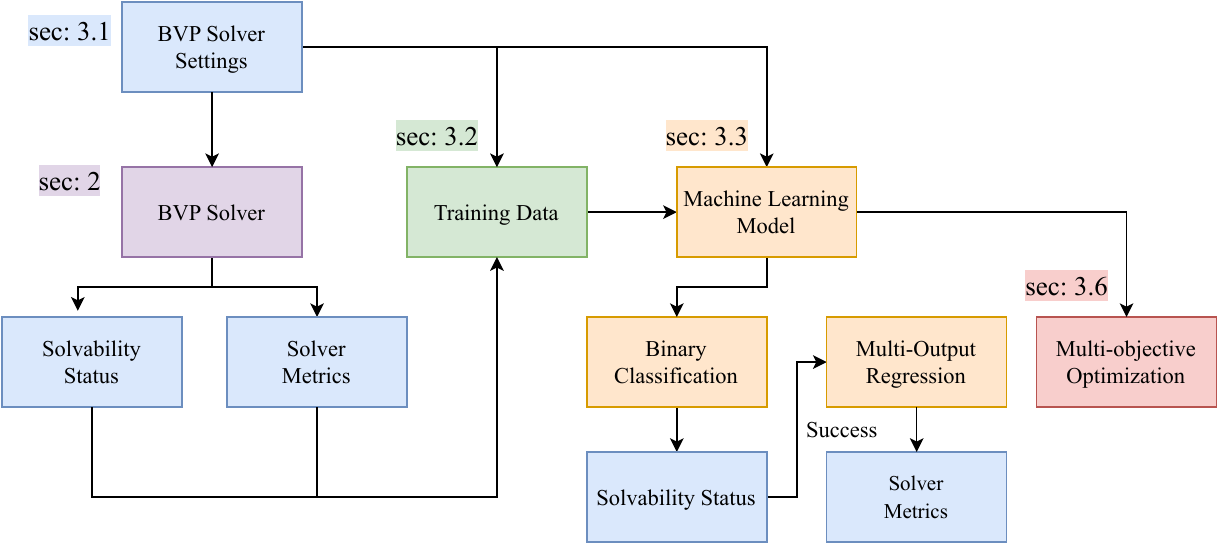}
    \caption{Workflow of machine learning-based optimization for tuning the numerical settings of BVP solvers. Different stages are colored corresponding to the section numbers in this chapter.}
    \label{fig:fine_tune}
\end{figure}

\subsection{Feature Selection}
Our feature selection process involved discussion with the domain experts to incorporate their experience in choosing the features based on the type of boundary value problems under inspection. The selected features are also confirmed using exploratory data analysis and the machine learning feature selection process. In the former, the correlations of the features within themselves and to the output features are considered, and in the latter, methods such as feature elimination and dimensionality reduction are utilized. We found a good agreement between the expert-based and ML-based feature selection processes.

\subsubsection{Input features}
We identified nine input features that affect the solvability and performance of the boundary value problem solvers. As discussed before, the BVP solver under consideration internally uses the Newton-Armijo method to solve the system of algebraic equations. Therefore, the input features consist of the generic BVP solver settings and the specific Newton-Armijo solver settings. The selected input features are as follows.

\begin{enumerate}

    \item { \emph{Test Case Type} is the type of the boundary value problem chosen from the test bench explained before. This is a categorical parameter where each category represents a problem from~\cite{soetaert2010package}.}
    
    \item { \emph{Newton Tolerance} handles the termination of the Newton solver. If the residual of the equation system is below the \emph{NewtonTolerance} the algorithm is considered converged.}
    
    \item { \emph{Newton Critical Tolerance} is an additional tolerance that is considered when the solver does not reach the more restrictive \emph{NewtonTolerance}. If the residual is below the \emph{NewtonCriticalTolerance} after the maximum number of iterations it is still considered as converged, otherwise as diverged.}

    \item { \emph{Newton Maximum Iterations} is the maximum number of iterations for the Newton solver.}
    
    \item { \emph{Newton Armijo Probes}  is the maximum Armijo probes used in the Newton-Armijo algorithm.}

    \item { \emph{Maximum Grid Points} is the maximum number of mesh points allowed when solving the boundary value problem.}

    \item { \emph{Add Factor} is the feature that controls the addition of new mesh points based on the difference between computed residual in the interval and the relative tolerance of the solver.}

     \item { \emph{Remove Factor} is the feature that controls the removal of mesh points if it is too fine based on the difference between the computed residual in the interval and the relative tolerance of the solver.}

    \item { \emph{Use Collocation Scaling} is the boolean variable that specifies whether the inner (non-boundary) conditions are scaled regarding the number of mesh points.}
\end{enumerate}

\subsubsection{Output features}
The output features represent the solvability status of the boundary value problem and the statistics that affect the time taken to solve the problem and the resources required. The four chosen output features are:
\begin{enumerate}

    \item { \emph{Success Status} is the boolean feature that specifies the solvability status of the BVP solver.}

    \item { \emph{Number of ODE Evaluations} is the number of times the solver calls the right-hand side function of the differential equation 
    in the algorithm.}

    \item { \emph{Number of Grid Points} is the number of grid points of the numerical solution.}

    \item { \emph{Maximum Residuum} is the maximum residual of the numerical solution.}
    
\end{enumerate}

\subsection{Data Collection}
As discussed before, a BVP solver based on MATLAB's bvp4c was used to collect the data for training the machine learning models. We set realistic ranges for the input features based on discussions with the domain experts.
Latin Hypercube Sampling was used as the sampling technique for creating the input data as this technique distributes samples evenly over sample space. Based on the boundary value problems under investigation and the domain expertise we chose either linear scale or logarithmic scale for the sampling. Table~\ref{tab:1} shows the ranges and scales used for each input feature for sampling. The generated input data are fed to the solver to generate the output labels. We generated 100,000 data points for each of the ten reference problems.

\begin{table*}[!t]
\caption{Ranges and scales used for sampling input features. The table also shows the default values of the features.}
\label{tab:1}
\begin{tabular}{|p{3.5cm}|p{3.5cm}|p{2.25cm}|p{2cm}|}
 \hline
 Input Feature & Feature Range & Default Value & Sampling Scale \\
 \hline
Maximum Grid Points & $[100, 10000]$  & \num{1000} & Linear \\
Newton Critical Tolerance & $[\num{e-12},100]$  & \num{1} & 
Logarithmic \\
Newton Armijo Probes & $[1, 10]$  & \num{2} & Linear \\
Newton Maximum Iterations & $[1, 100]$  & \num{4} & Linear \\
Newton Tolerance & $[\num{e-12}, \num{e-2}]$  & \num{e-12} & Logarithmic \\
Add Factor & $[1, 1000]$  & \num{100} & Logarithmic \\
Remove Factor & $[0, 2]$  & \num{0.5} & Linear \\
Use Collocation Scaling & \emph{\{True, False\}}  & \emph{True} & Discrete \\
Test Case Type & Ten discrete problem types   &---& Discrete \\
 \hline
\end{tabular}
\end{table*}

\subsection{Model Selection}
Based on our problem statement, we evaluated several binary classification and multi-output regression algorithms. The selected algorithms are classified into five major categories which are support vector machine, neighbors-based, boosting, bagging, and artificial neural net-based algorithms. One best algorithm from each category is selected based on the accuracy of the test dataset. The description of the chosen categories is listed as follows.

\begin{enumerate}

    \item {Support Vector Machine (SVM)~\cite{noble2006support} is a supervised learning algorithm that finds the best hyperplane to separate the features into different domains. They solve the non-linear problems with the help of kernels which implicitly map the input features into a high dimensional space.}

    \item {Boosting algorithms are based on the gradient boosting technique, in which the classification model is constructed as an ensemble of weak models. It is a sequential process, where each subsequent model attempts to correct the errors of the previous model. We selected the light gradient boosting machine (LGBM) algorithm~\cite{ke2017lightgbm} that includes sparse optimization, parallel training, multiple loss functions, regularization, and early stopping techniques to provide an accurate and generalized model.}

    \item {Bagging algorithms are based on ensemble techniques, where the dataset is randomly sampled into subsets and each of the subsets is trained using different models. The predictions are calculated by combining individual subset predictions.
    Random Forest (RF)~\cite{breiman2001random} is selected in this category which is a bagging or bootstrap aggregation technique based on decision trees.}

    \item {Neighbors-based algorithms work on the fundamental idea of finding the $n$ training samples closest in distance to a new point and assigning the label from these. Euclidean distance is the most common metric used to measure distance. These are non-parametric algorithms that do not make strong assumptions about the form of the mapping function. Hence, they frequently work in scenarios where the decision border is not regular. We selected the K-nearest neighbor algorithm ~\cite{fix1989discriminatory},~\cite{cover1967nearest},  that works based on the $K$ nearest neighbors of each data point.}

    \item {Artificial neural network (ANN) is a fully connected feed-forward network that makes use of an aggregation of functions to understand and convert input data that is presented in one form into the desired output. A traditional neural network's fundamental building blocks are nodes, which are arranged into layers. To get the final forecast, the input characteristics are transferred through these layers using a series of non-linear operations. For binary classification, we used a multilayer perceptron (MLP) of two hidden layers with 100 neurons each. For regression, we used an MLP comprising six hidden layers with 250 neurons each.}
    
\end{enumerate}

\subsection{Training and Testing}
We divided the dataset into a \SI{80}{\percent}  training data and a \SI{20}{\percent} testing set. The training data is further divided into a training set and a validation set. The training set is used to train the regression models and the validation set is used to tune the model hyper-parameters. The testing set is used for unbiased evaluation of the model. As a part of data preparation, we performed input feature scaling and output feature transformation to tailor the data for the machine learning models. During the exploratory data analysis, we observed that the output data had a skewed distribution with chunks of repeated values and outliers. Hence we transformed the distribution to follow a normal distribution that spread out the most frequent values and also reduces the impact of outliers.

\subsection{Evaluation Metrics}
In this section, we list the metrics that are used to evaluate the classification and regression models. The metrics are selected based on the type of algorithms, scale, and distribution of the input and output features that reveal the fairness, robustness, and reliability of ML models.

\subsubsection{Binary Classification}
For the equations used below, we define TP, FP, TN, and FN as true positives, false positives, true negatives, and false negatives respectively.

\begin{enumerate}

    \item {Accuracy: It measures how many observations were correctly classified overall. This is expressed as a percentage. Accuracy is calculated as below.}

    \begin{align*}
    P &= \frac{TP + TN}{TP + FP + TN + FN}
    \end{align*}

    \item {Precision (P): It measures the proportion of positive classifications by the model that were actually correct. Precision is calculated as below.}

    \begin{align*}
    P &= \frac{TP}{TP + FP}
    \end{align*}

    \item {Recall (R): It measures the proportion of actual positives classifications that were identified correctly by the model. Recall is calculated as below.}

    \begin{align*}
    R &= \frac{TP}{TP + FN}
    \end{align*}
    
\end{enumerate}

\subsubsection{Multi-output Regression}
For the equations used below, we define $n$ as the number of data points, $y_i$, and $\hat y_i$ as the actual and predicted values respectively for the data point $i$ and $\bar y$ as the mean value of the data points.

\begin{enumerate}

\item {Mean absolute percentage error (MAPE) is a statistical measure to
    evaluate the accuracy of a regression model. The error is independent of the scale of the output as it measures the accuracy as a percentage. MAPE is calculated as below.}

    \begin{align*}
    MAPE &= \frac{100}{n} \sum_{i=1}^{n} \left| \frac{y_i-\hat y_i}{y_i}\right| \\
    \end{align*}
    
\item {Root Mean squared error (RMSE) measures the square root of the average squared errors. The RMSE is a good estimate for ensuring that the ML model has no outlier predictions with large errors since it puts a higher weight on these errors due to the squaring. RMSE is calculated as below.}
    
    \begin{align*}
    RMSE &= \sqrt{\frac{\sum_{i=1}^{n} (y_i-\hat y_i)^2}{n}}
    \end{align*}

    \item {Coefficient of determination ($R^2$ Score) is the measure of how
    close the data points are to the fitted regression line. It explains how much of the variance of actual data is explained by the predicted values. $R^2$ Score is calculated as below.}

    \begin{align*}
    R^2(y, \hat y) &= 1 -  \frac{\sum_{i=1}^{n} (y_i - \hat y_i)^2}{\sum_{i=1}^{n} (y_i - \bar y_i)^2}
    \end{align*}
\end{enumerate}

\begin{table}[!t]
\caption{Functions provided by the bvpTune library}
\label{tab:4}
\begin{tabular}{|p{6.3cm}|p{5cm}|}
 \hline
 Function & Description \\
 \hline
 getSolvabiltyStatus(TestCaseType, Setting)
 &  returns whether the numerical setting is solvable or not for the specified TestCaseType of the boundary value problem. \\
 \hline
  getSolverPerformance(TestCaseType, Settings)
 &  returns the predicted solver statistics of the numerical setting for the specified TestCaseType of the boundary value problem. \\
 \hline
 \begin{tabular}{l}getOptimalODEevals(TestCaseType) \\ getOptimalGridPoints(TestCaseType) \\
 getOptimalResiduum(TestCaseType) \end{tabular} & returns the optimal number of ODE evaluations, grid points, and residuum for the specified TestCaseType of the boundary value problem.  \\
 \hline
 \begin{tabular}{l}getOptimalODEevalsAndGridPoints(TestCaseType) \\ getOptimalGridPointsAndResiduum(TestCaseType) \\
 getOptimalResiduumAndODEevals(TestCaseType) \end{tabular}
 & return pareto front numerical settings as a data frame that optimizes the two criteria for the specified TestCaseType of the boundary value problem.  \\
 \hline
 getOptimizedSettings(TestCaseType)
 & returns pareto front numerical settings as a data frame that optimizes the ODE evaluations, grid points, and the residuum
  for the specified TestCaseType of the boundary value problem.  \\
 \hline
visualize(dataframe)
 & visualize the 2d/3d pareto front plots for the provided data frame returned by the optimization functions.  \\
 \hline
\end{tabular}
\end{table}

\subsection{Multi-objective Optimization}
\label{sec:4}

After determining the best machine learning models based on the evaluation metrics, we use them for optimization. The objective in our problem setting is to minimize the residual error obtained from the BVP solver by optimizing the computational resources such as the number of ODE evaluations and the grid points. Hence, we define the problem as a multi-objective optimization as shown in equation~\eqref{eq:1}.

    \begin{equation}
    \min_{ x \in X} (\phi_{eval}(x), \phi_{gp}(x), \phi_{res}(x))
    \label{eq:1}
    \end{equation}
    subject to the constraints,
    \begin{align*}
    f_{lb}(i) \leq x_i \leq f_{ub}(i), \mbox{ for } i=1,...,n 
    \end{align*}
    and
    \begin{align*}
    \phi_{cl}(x) = 1
    \end{align*}
where $\phi_{eval}$, $\phi_{gp}$, and $\phi_{res}$ are the selected ML regression algorithms for predicting ODE evaluations, grid points, and maximum residuum respectively. $\phi_{cl}$ is the selected ML classification algorithm.
$X$ is the feasible set of input vectors and
$x = [x_1, x_2,....,x_n]$. The functions $f_{lb}(i)$ and $f_{ub}(i)$ provide the upper and lower bounds of the $i^{th}$ input parameter search space respectively. Along with equation \eqref{eq:1} which optimizes all three objectives, we also examine one and two objective optimizations to cater to the requirements based on the different problem setups where one or two objectives can be relaxed.

We use the Optuna library~\cite{optuna_2019} for solving our multi-objective optimization. This library provides a framework to automatically search for the optimal input vectors and efficiently search large spaces and prune unpromising results using state-of-the-art algorithms such as Nondominated Sorting Genetic Algorithm II~\cite{deb2002fast}, the Tree-structured Parzen Estimator Algorithm~\cite{ozaki2020multiobjective, ozaki2022multiobjective}, and the Quasi Monte Carlo Sampling Algorithm~\cite{bergstra2012random}, etc. We provide the objective function defined in the equation \eqref{eq:1} and the best machine learning model to the Optuna framework for optimization.

\begin{figure}
    \centering
    \includegraphics[width=\linewidth]{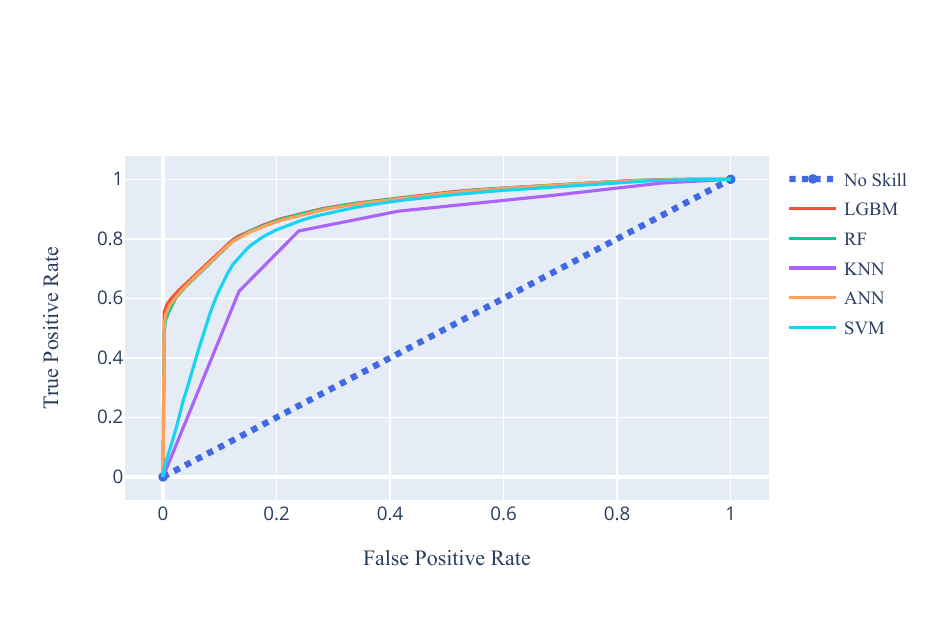}
    \caption{ROC curves for different classification algorithms}
    \label{fig:roc_curve}
\end{figure}

\section{Experiments and Results}
\label{sec:5}
In this section, we present the performance of the proposed workflow on the BVP solver. The results reveal the necessity and effectiveness of our workflow in fine-tuning the numerical settings of the BVP solvers. We further discuss the different quality aspects related to the proposed workflow to show its scalability, reliability, and stability. The scalability aspects are computed on a workstation with a 40-core Intel® Xeon® E5-2680 v2 (2.80 GHz) CPU.

\begin{table*}[!t]
\caption{Performance of classification models on the dataset}
\label{tab:2}
\begin{tabular}{|p{4.5cm}|p{2.5cm}|p{2.25cm}|p{2cm}|}
 \hline
 ML Algorithm & Accuracy & Precision & Recall \\
 \hline
Support Vector Machine & 0.8458  & 0.8722 & \textbf{0.9355} \\
Light Gradient Boosting Machine  & \textbf{0.8593} & \textbf{0.9048} & 0.9120 \\
Random Forest & 0.8123  & 0.8696 & 0.8874 \\
K Nearest Neighbor & 0.8206  & 0.875 & 0.893 \\
Artificial Neural Net & 0.8562  & 0.9023 & 0.9106 \\
 \hline
\end{tabular}
\end{table*}

\subsection{Binary Classification}
\begin{figure}
    \centering
    \includegraphics[width=\linewidth]{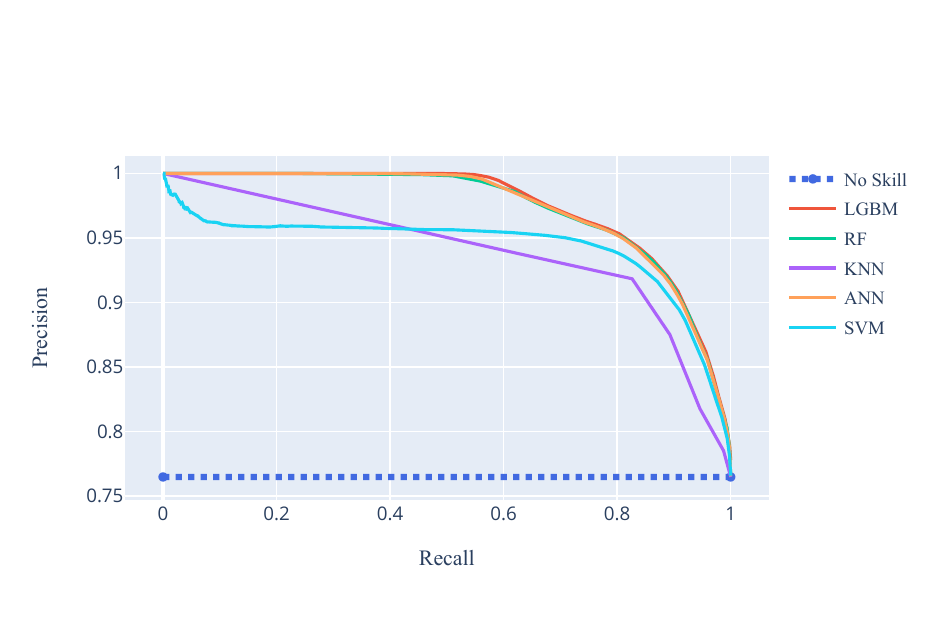}
    \caption{Precision-Recall curves for different classification algorithms}
    \label{fig:pr_curve}
\end{figure}
The machine learning algorithms to achieve binary classification were trained independently of the types of boundary value problems discussed in Section~\ref{sec:2}. We were able to generalize the models for all the considered types of problems given to the BVP solver. Table~\ref{tab:2} shows the performance of the binary classification models on the test dataset. While the LGBM algorithm provided the best results in terms of accuracy and precision, the SVM algorithm provided the best recall. For further investigation, we plotted ROC curves of the ML models, as seen in Figure~\ref{fig:roc_curve}. These curves show the performance of classification models at different classification thresholds. The area under the ROC curve (AUC) measures how well predictions are ranked, rather than their absolute values. As observed from the figure, LGBM and ANN have better AUCs.

\begin{table}[!t]
\caption{Performance of regression models on the dataset}
\label{tab:3}
\begin{tabular}{|p{4.5cm}|p{2.5cm}|p{2.25cm}|p{2cm}|}
 \hline
 ML Algorithm & RMSE & MAPE & $R^2$ Score \\
 \hline
Support Vector Machine & 10658.80  & 229.40 & 0.094 \\
Light Gradient Boosting Machine  & 2760.38 & 4.30 & 0.95 \\
Random Forest & \textbf{1870.41}  & \textbf{1.31} & \textbf{0.97} \\
K Nearest Neighbor & 8315.08  & 69.65 & 0.48 \\
Artificial Neural Net & 3954.85  & 7.90 & 0.90 \\
 \hline
\end{tabular}
\end{table}

The dataset generated by the BVP solver consisted of more positive cases (solvable settings) compared to the negative cases (non-solvable settings). To evaluate the classification models with this imbalanced dataset, we plotted precision-recall curves as depicted in Figure~\ref{fig:pr_curve}. Also like ROC curves, precision-recall curves provide a graphical representation of a classifier’s performance across many thresholds rather than a single value. The blue dotted line represents a baseline classifier that simply predicts that all instances belong to the positive class. As we can see from the figure, LGBM, RF, and ANN have the largest area under the precision-recall curve (AUC-PR).

\begin{figure}
    \includegraphics[width=\linewidth]{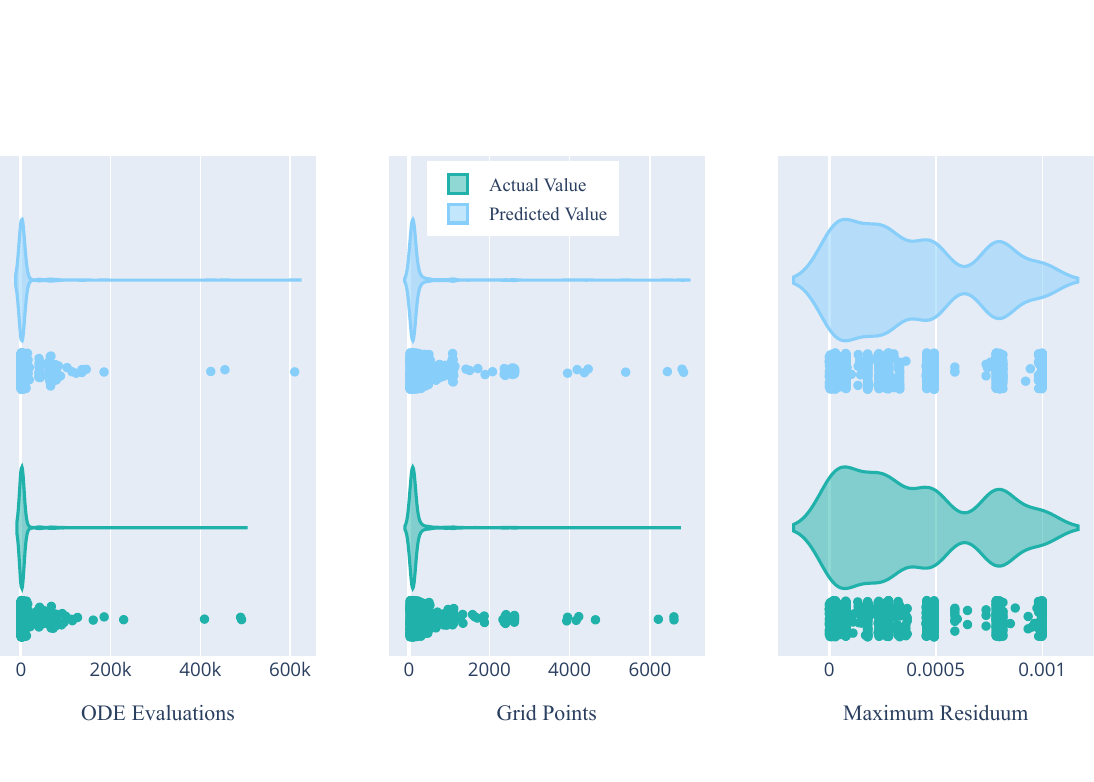}
    \centering
    \caption{Violin plot comparing the distribution of the actual and predicted values of ODE evaluations, grid points and maximum residuum.}
    \label{fig:violin_plot}
\end{figure}

\subsection{Multi-output Regression}
As discussed in the proposed workflow, we train regression models to predict the solver's performance after estimating the solvability status of the boundary value problems using a binary classification technique. During the exploratory data analysis, we discovered that the solver-specific numerical settings have little (problem 7) to no effect (problems 1, 3, and 4) on the performance of linear problems. The effect is significant for non-linear problems. In contrast to classification, regression techniques strongly depend on the problem type to forecast its performance. Hence, the regression models are trained based on the type of problem that is solved. Table~\ref{tab:3} shows the problem-specific performance of the regression models. The table shows that the RF algorithm outperformed other algorithms in terms of evaluation metrics. The comparison of the output distribution of the actual and predicted values using the RF algorithm is shown in Figure~\ref{fig:violin_plot}.

\begin{figure}
    \centering
    \includegraphics[width=\linewidth]{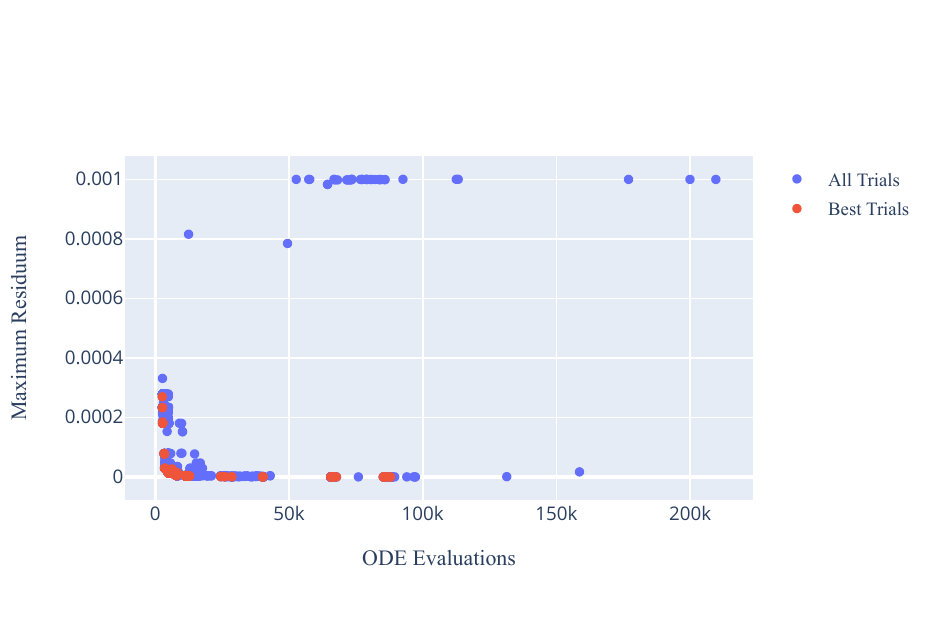}
    \caption{Pareto Front plot showing the results of two-objective optimization trials along the number of ODE evaluations and the maximum residuum.}
    \label{fig:2_opt}
\end{figure}

\subsection{Multi-objective Optimization}
\begin{figure}
    \centering
    \includegraphics[width=\linewidth]{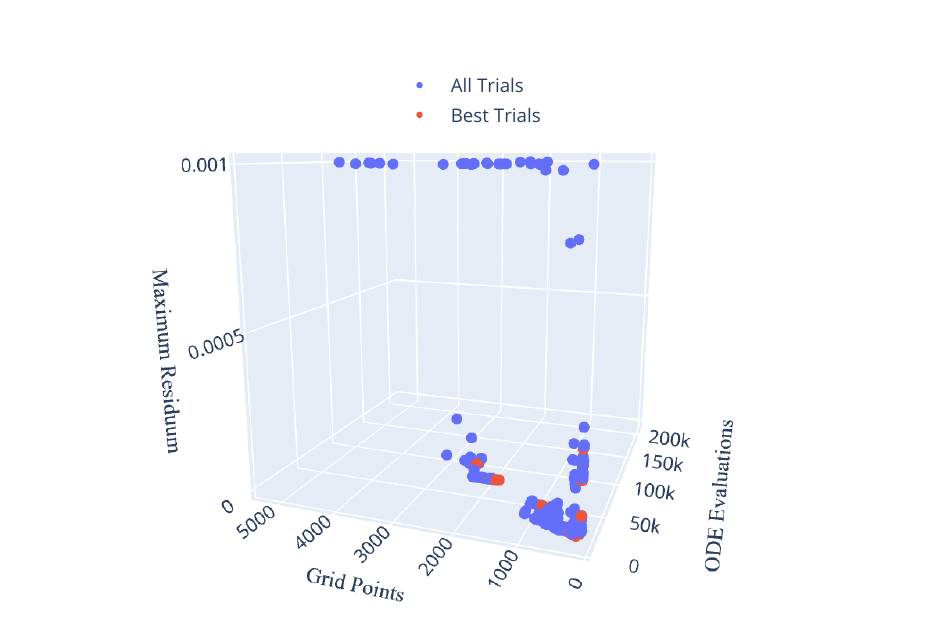}
    \caption{Pareto Front plot showing the results of three-objective optimization trials.}
    \label{fig:3_opt}
\end{figure}

The ML models, which are chosen based on their performance, are used to optimize the problem~\eqref{eq:1}. The equation can be customized based on the requirements to perform single-, two-, and three-objective optimizations. Figure~\ref{fig:2_opt} shows the Pareto front for optimizing the number of ODE evaluations and the maximum residuum. The 3D Pareto front for all three criteria optimization is shown in 
Figure~\ref{fig:3_opt}. In these plots, the blue dots represent the regular optimization trials, and the red dots represent the best trials chosen. Since predictions from the ML models are calculated in real-time, the time required for optimization is reduced. To examine the reliability of the results produced from our workflow, the best trials of the optimization were tested against the actual output from the BVP solver.

\subsubsection{Fine-tuning library for the BVP solver}
Based on the proposed workflow in the paper, we implemented a Python 
library, bvpTune, to fine-tune the numerical settings of the BVP solver. Using the library, users can retrieve the solvability of the desired problem based on the solver settings. The library also facilitates the optimization of the quality metrics based on the custom queries provided through the function interfaces. The library can be installed as a Python package. The source code and documentation are provided on GitHub.\footnote{\url{https://github.com/VictorVinySaajan/bvpTune}} The functionalities provided by the library are described in Table~\ref{tab:4}.

\section{Workflow Evaluation}
\begin{figure}
    \includegraphics[width=\linewidth]{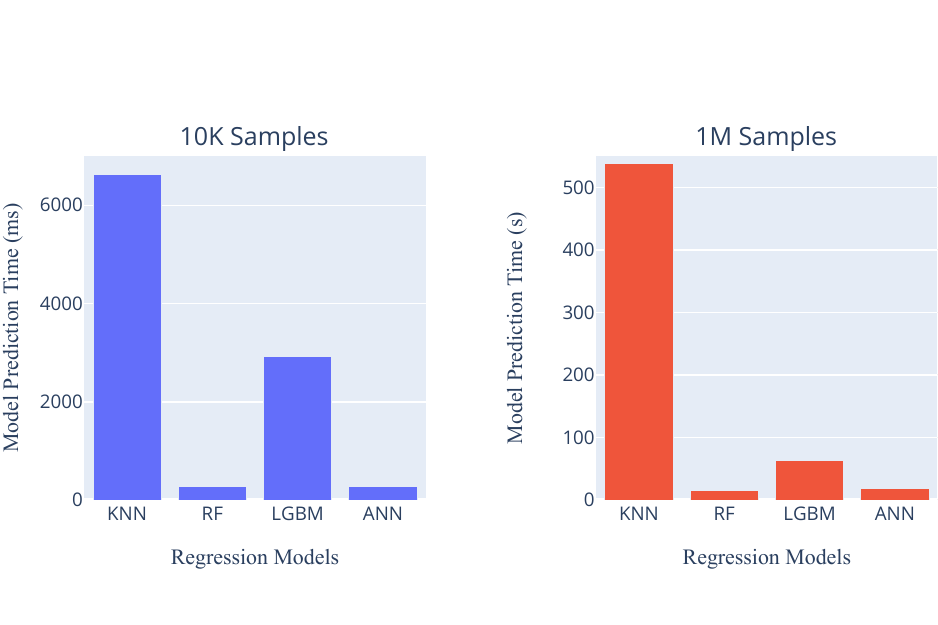}
    \centering
    \caption{Bar plot showing the scalability of the regression models.}
    \label{fig:time_rgr}
\end{figure}

In this section, we evaluate the scalability, stability, and reliability of our workflow for fine-tuning the numerical settings of the BVP solvers. The scalability of any workflow is important for real-world applications. The scalability of our workflow mainly depends on the scalability of the ML model. Hence, in our context, "scalability" refers to the ML model's capacity for handling large volumes of data and performing numerous computations rapidly and effectively. Since the training is only performed once, we primarily concentrate on the prediction time. We evaluated the computing time of the predictions for a significant data set in order to ascertain whether the models scale well with the data. 
Figures~\ref{fig:time_clf} and ~\ref{fig:time_rgr} compare the time taken by the classification and regression models for predicting 10,000 and 1,000,000 outputs, respectively. We did not include the SVM algorithms in the analysis as they did not scale well with the data. The figures show that, with the exception of the KNN, all the models scale reasonably well with the data for both classification and regression. We observe that the LGBM algorithm is the fastest for classification and the RF algorithm is the fastest for regression. To examine the reliability of the workflow, the outputs from the bvpTune library are tested against the actual outputs from the BVP solver. We found that both observations were in agreement with each other.

As mentioned previously, the regression models are trained specifically to the type of reference boundary value problems. We do this because the effect of the numerical settings on the solver's performance varies significantly as the problem changes. To confirm this, we compared the output distributions of our reference problems for the same set of inputs. We quantified the drift using the population stability index ($PSI$)~\cite{yurdakul2018statistical} that measures population stability between two distribution samples. For samples divided into $B$ bins, the $PSI$ is calculated as below:

\begin{figure}
    \includegraphics[width=\linewidth]{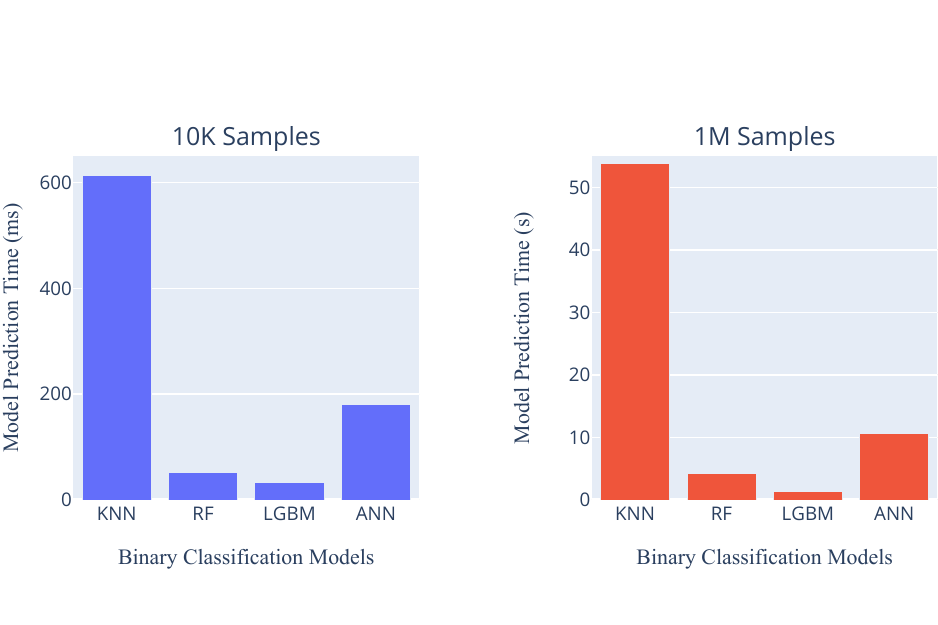}
    \centering
    \caption{Bar plot showing the scalability of the classification models.}
    \label{fig:time_clf}
\end{figure}

\begin{equation}
    PSI = \sum_{i=1}^{B} (y_i - y_{ri}) \ln(y_i/y_{ri})
\end{equation}

\begin{table*}[!t]
\caption{Pairwise combination of PSI values for the reference boundary value problems defined in~\cite{soetaert2010package}. The values of $i^{\text{th}}$ row and $j^{\text{th}}$ column are calculated by keeping output distributions of Reference Problem Number(i) and Reference Problem Number(j) as initial populations respectively.}
\label{tab:5}
\centering
\begin{tabular}{|M{1.4cm}|M{1.3cm}|M{1.3cm}|M{1.3cm}|M{1.3cm}|M{1.3cm}|M{1.3cm}|M{1.3cm}|}
 \hline
 Reference Problem Number & 7 & 19 & 20 & 22 & 23 & 24 & 33 \\
 \hline
7 & 0.0000  & 2.5384 & 4.5560 & 2.4123 & 1.5956 & 2.5570 & 1.6833 \\
19 & 6.0619 &  0.0000 &	1.9525 & 4.2577	& 3.9291 & 6.2365 &	2.0148 \\
20 & 5.2940 &	0.6268 & 0.0000 & 4.6000 & 2.8879 &	6.6869 & 1.4494 \\
22 & 5.7805 &	3.8481 &  3.5217 & 0.0000	& 2.1835 & 1.3254 &	3.0429 \\
23 & 3.3124 &	2.2905 & 2.5908 & 2.1386 & 0.0000 &	0.7574 & 3.6245 \\
24 & 4.2097 & 3.6440 & 3.6107 &	1.5194 & 0.7938 & 0.0000 & 4.3400 \\
33 & 5.7000 &	4.9902 & 5.8062 & 6.6005 & 7.5409 &	7.0110 & 0.0000 \\
 \hline
\end{tabular}
\end{table*}

The $y_i$ is the proportion of the new output distribution that falls in the $i^{th}$ bin and $y_{ri}$ is the proportion of the initial output distribution that falls in the $i^{th}$ bin. According to the widely used rule of thumb proposed in~\cite{siddiqi2012credit}, the drift between the populations is stable if the $PSI$ value is less than 0.1. A value between 0.1 to 0.25 indicates a small shift between the populations and a higher than 0.25 value represents two substantially different populations. The PSI values observed from Table~\ref{tab:5}, suggest that all the reference problems have significantly different behaviors corresponding to the same inputs and cannot be generalized.

\section{Conclusion}
In this chapter, an ML-based optimization workflow for tuning the numerical settings of boundary value problem solvers is proposed, and models based on binary classification and multi-output regression are established. We demonstrated the optimization of quality metrics exemplarily for a BVP solver using these models. Furthermore, several machine learning algorithms for parameter tuning are explored, and the scalability and reliability of these models are discussed. Experimental results show that
LGBM has good accuracy and computational performance for classification, while
Random Forests have good accuracy and computational performance for regression. 
Additionally, the stability of the models for unseen problem types is discussed and can be continued for different boundary value problem solvers in future work.
\label{sec:6}

\section{Acknowledgement}
This work was developed in Transport Processes department of Fraunhofer ITWM. Furthermore, we acknowledge the support of the High Performance Center for Simulation and Software-Based Innovation.
\label{sec:7}

\section{Authors Contribution}
\label{sec:8}

The authors confirm their contribution to the paper as follows: study conception and design:  Viny Saajan Victor, Manuel Ettm\"{u}ller, Andre Schmei{\ss}er, Heike Leitte and Simone Gramsch; data collection: Viny Saajan Victor and Manuel Ettm\"{u}ller; analysis and interpretation of results: Viny Saajan Victor; draft manuscript preparation: Viny Saajan Victor. All authors reviewed the results and approved the final version of the manuscript. 

\section{Conflict of Interest}
All the named authors of this chapter declare that they have no conflicts of interest, financial or otherwise.
\label{sec:9}

%
%
%
%


%

\backmatter


\end{document}